\documentclass[a4paper,fleqn]{cas-sc}

\usepackage[authoryear,longnamesfirst]{natbib}

\newtheorem{theorem}{Theorem}[section]
   
\newtheorem{remark}[theorem]{Remark}

\newtheorem{lemma}[theorem]{Lemma}
\newdefinition{rmk}{Remark}
\newproof{pf}{Proof}
\newproof{pot}{Proof of Theorem \ref{thm2}}
\newcommand{\E}{\mathbf{E}}
\renewcommand{\P}{\mathbf{P}}
\newcommand{\Prob}[1]{\P\left\{#1\right\}}

\newcommand{\N}{\mathbb{N}}
\newcommand{\R}{\mathbb{R}}

\newcommand{\sB}{\mathcal{B}}
\newcommand{\sP}{\mathcal{P}}
\newcommand{\sR}{\mathcal{R}}
\newcommand{\sM}{\mathcal{M}}

\newcommand{\one}{\mathbf{1}}

\newcommand{\eps}{\varepsilon}

\newcommand{\ns}{\lfloor ns \rfloor}
\newcommand{\nz}{\lfloor nz \rfloor}
\newcommand{\dGH}{\mathsf{d}_{\mathrm{GH}}}
\newcommand{\dH}{\mathsf{d}_{\mathrm{H}}}
\newcommand{\nt}[1][t]{\lfloor n#1\rfloor}

\newcommand{\salg}{\mathfrak{F}}

\usepackage{fancybox}
\newlength{\querylen}
\begin{document}
\let\WriteBookmarks\relax
\def\floatpagepagefraction{1}
\def\textpagefraction{.001}
\shorttitle{Random Bridges in Spaces of Growing Dimension}
\shortauthors{Bochen Jin}

\title [mode = title]{Random Bridges in Spaces of Growing Dimension}

\author[1]{Bochen Jin}[type=editor,
                       orcid=0009-0008-3394-2973]
\ead[url]{bochen.jin@unibe.ch}

\affiliation[1]{organization={Institute of Mathematical Statistics
    and Actuarial Science, University of Bern},
                addressline={Alpeneggstrasse 22}, 
                postcode={3012 Bern}, 
                country={Switzerland}}

\begin{abstract}
  We investigate the limiting behaviour of the path of random
  bridges treated as random sets in $\mathbb{R}^{d}$ with the
  Euclidean metric and the dimension $d$ increasing to infinity. The
  main result states that, in the square integrable case, the limit (in
  the Gromov-Hausdorff sense) is deterministic, namely, it is $[0,1]$
  equipped with the pseudo-metric $\sqrt{|t-s|(1-|t-s|)}$. We also
  show that, in the heavy-tailed case with summands regularly varying
  of order $\alpha \in (0,1)$, the limiting metric space has a random
  metric derived from the bridge variant of a subordinator.
\end{abstract}


\maketitle

\section{Introduction}
Consider a random walk 
\begin{displaymath}
  S_n^{(d)} := X_1^{(d)} + X_2^{(d)} + \cdots + X_n^{(d)}, \quad n \geq 1,
\end{displaymath}
where $X_i^{(d)} := (X_{i,1}^{(d)}, \ldots, X_{i,d}^{(d)})$, $i\geq1$,
are i.i.d.\ random vectors in $\R^d$.  Set $S_0^{(d)}=0$.  The path of
this random walk is a finite random subset of $\R^d$, which is treated
as a finite metric space
\begin{displaymath}
  \sM_n^{(d)}:=\Big\{ S_0^{(d)},\ldots,S_n^{(d)}\Big\}
\end{displaymath}
with the $\ell^2$ metric induced from $\R^d$.

If the dimension $d$ is fixed and $n\to\infty$, then the random set
$n^{-1/2}\sM_n^{(d)}$ converges in distribution to the path of the
Brownian motion in $\R^d$ under the assumptions $\E
X_1^{(d)}=0$ and $\E \|X_1^{(d)}\|^2=1$. This convergence is
understood as convergence of random compact sets in distribution in
the topology generated by the Hausdorff distance.

In order to study the convergence when the dimension varies, it is
possible to use the concept of the Gromov--Hausdorff distance between
metric spaces. The \textit{Gromov--Hausdorff distance} between metric
spaces $\mathbb{X}:=(X, \rho_X)$ and $\mathbb{Y}:=(Y,
  \rho_Y)$ is defined as
\begin{displaymath}
  d_{GH}(\mathbb{X}, \mathbb{Y}) := \inf_{i: X \hookrightarrow Z, j: Y \hookrightarrow Z} d_H
  \big(i(X), j(Y) \big),
\end{displaymath}
where the infimum is taken over all isometric embeddings $i$ and $j$
into all possible metric spaces $(Z, d)$ which can embed $X \
\text{and} \ Y$. The Hausdorff distance between sets $A$ and $B$
in $(Z,d)$ is defined as
\begin{displaymath}
  d_H(A, B) := \inf \{ \epsilon > 0: A \subset B^{\epsilon} \
  \text{and}\ B \subset A^{\epsilon} \},
\end{displaymath}
where $A^{\epsilon} := \{x: d(x, A) < \epsilon \}$ is the $\epsilon$-neighbourhood of $A$,
see \citep[Chapter 7]{Buragoetal2001}. 
This is indeed a metric on the family of all compact metric
spaces (considered equivalent up to an isometry), which generates the
Borel $\sigma$-algebra and so makes it possible to define a random
metric space.

\cite{KabluchkoMarynych2024} showed that, assuming
finiteness of $\E\|X^{(d)}\|^2$, the random metric space 
$\big(n^{-1/2}\mathcal{M}_n^{(d)}, \|\cdot\|\big)$ converges in probability to a
deterministic limit in the Gromov-Hausdorff sense.
The limiting space $W$ is the so-called \textit{Wiener
spiral}, which is the set $\big\{ \one_{[0, t]} : 0 \leq t \leq 1 \big\}$ 
of indicator functions on $[0, 1]$, considered as a subset of the 
Hilbert space $L^2[0,1]$ and endowed with the $\ell^2$ metric. 
The Wiener spiral is isometric to $[0,1]$ with the metric $r(t,s)=\sqrt{|t-s|}$.
The key argument for this result is the uniform law of large numbers
which says that the squared distances
$n^{-1}\|S_{\nt}^{(d)}-S_{\ns}^{(d)} \|^2$ converge uniformly to $|t - s|$ over
$t, s \in [0, 1]$, and so the non-squared distance becomes $\sqrt{|t-s|}$ in the
limit. In the subsequent paper \cite{Kabluchkoetal2024}, 
the case of non-integrable increments was considered. 
In this case the limit is a random metric space.

It is possible to consider variants of the basic random walk obtained by applying linear
transformations. The main example is the bridge variant given by
\begin{displaymath}
  B_k^{(d)}:=S_k^{(d)}-\frac{k}{n}S_n^{(d)}, \quad
  1\leq k\leq n,\, n\in\N.
\end{displaymath}

Let $(\sB_n^{(d)},\|\cdot\|)$  be a random metric space derived from the path of
the bridge variant random walk 
\begin{displaymath}
  \sB_n^{(d)}:=\Big\{B_0^{(d)}, \ldots,B_n^{(d)}\Big\}
\end{displaymath}
and  endowed with the Euclidean metric.
We study the limit of $\sB_n^{(d)}$.
In Section 2, we deal with the square integrable case. In Section 3,
we address the case of the infinite mean. Proofs of the theorems in
Section 2 and Section 3 are presented in Section 4 and Section 5,
respectively.

\section{Convergence to the bridge variation of the Wiener spiral}
We adapt the following conditions, which are imposed by 
\cite{KabluchkoMarynych2024} for studying the
Gromov-Hausdorff convergence of the random metric space 
$\big(n^{-1/2}\sM_n^{(d)},\|\cdot\|\big)$.
\begin{enumerate}[(a)]
\item The increment $X_1^{(d)}$ is centred, square integrable and
  normalized, that is, $\E X_1^{(d)} = 0\ \text{and}\ \E \| X_1^{(d)}\|^2=1$.
\item The components of $X_1^{(d)}$ are mutually uncorrelated, that is,
$\E \big[X_{1,j}^{(d)}X_{1,k}^{(d)} \big] = 0$ for all $j,k \in {1, \ldots, d}, \; j \neq k.$
\item The sequence $(\| X_1^{(d)} \|^2)_{d \in \N}$ is uniformly
  integrable, that is,
  \begin{displaymath}
    \lim_{A \to \infty} \sup_{d \in \N}
    \E \big(\| X_1^{(d)} \|^2 \one_{{\| X_1^{(d)} \|^2 > A }} \big)= 0.
  \end{displaymath}
\item The individual components of $X_1^{(d)}$ are negligible, that is,
  \begin{displaymath}
    \lim_{d \to \infty} \max_{k \in {1, \ldots, d}} \E\big(X_{1,k}^{(d)}\big)^2 = 0.
  \end{displaymath}
\end{enumerate} 

We now formulate our result for the bridge random walk. While the 
limit is a pseudo-metric space, it is still possible to use the 
Gromov-Hausdorff distance to formulate the convergence statement.

\begin{theorem}
  \label{the_br_l2}
  Let $n = n(d)$ be an arbitrary sequence of positive integers such
  that $n(d) \to \infty$ as $d \to \infty$. Suppose that all
  conditions (a)-(d) are fulfilled. Then, as $d \to \infty$, the
  random metric space $\big(n^{-1/2}\sB_n^{(d)},\|\cdot\|\big)$ converges in
  probability to the space $[0,1]$ with the pseudo-metric given by
  $\sqrt{(1 -|t-s|)|t-s|}$.
\end{theorem}

Let $(w_t)_{t\in[0,1]}$ be the standard Brownian motion.
The space $\big([0,1],\sqrt{(1 -|t-s|)|t-s|}\big)$ is isometric to the set of
random variables determined by the standard Brownian bridge
$\big\{b_t=w_t-tw_1:t\in[0,1]\big\}$, and so is indeed a pseudo-metric space.
The limit is also isometric to the bridge version of the Wiener
spiral, which is defined as the set
$\big\{ \one_{[0, t]} - t: 0 \leq t \leq 1 \big\}$ and
considered as a subset of $L^2[0, 1]$.

\begin{remark}
  \label{rem:1}
  The convergence of the finite point set is equivalent to the 
  convergence of the curves obtained by linear interpolation of these points. 
  In view of this, the validity of Theorem \ref{the_br_l2} can be alternatively confirmed
  by the continuous mapping theorem, with the continuous transformation
  \begin{equation}
    \label{eq:add_1}
    (\phi(t))_{t\in[0,1]}\mapsto(\phi(t)-t(\phi(1)))_{t\in[0,1]}
  \end{equation}
  applied to both linearly interpolated sequence of sums and the limiting Wiener spiral.
\end{remark}

\section{The case of non-integrable increments}
Suppose that the following conditions imposed by \cite{Kabluchkoetal2024} hold.
\begin{enumerate}[(I)]
\item There exist constants $(a(n))_{n \in \N}$ and a L\'evy measure
  $\nu$ on $(0, \infty)$ satisfying
  \begin{displaymath}
    \int_{(0, \infty)} \min(1,x) \nu(dx) < \infty, 
  \end{displaymath}
  and such that
  \begin{equation}
    \label{equ_ass_i_1}
    n \Prob{ (a(n))^{-1}\|X_1^{(d)}\|^2 \in\cdot}
    \overset{v}\to\nu(\cdot) \quad\text{as} \ d \to \infty,
  \end{equation}
  where $\overset{v} \to$ stands for the vague convergence of measures
  on $(0, \infty)$. Suppose further that
  \begin{equation}
    \label{equ_ass_i_2}
    \lim_{s \to 0+} \limsup_{d \to \infty} \frac{n}{a(n)} \E
    \Big(\| X_1^{(d)} \|^2 \one_{\| X_1^{(d)} \|^2 \leq sa(n)}\Big)=0.
  \end{equation}
\item For all fixed $s > 0$ and $\eps > 0$,
  \begin{equation}
    \label{equ_ass_ii}
    \lim_{d \to \infty} \Prob{ \Big|\big
      \langle\Theta_1^{(d)},\Theta_2^{(d)} \big\rangle\Big|>\eps
      \Bigm| \| X_1^{(d)} \|^2 \geq s a(n), \| X_2^{(d)} \|^2 \geq sa(n)}=0,
  \end{equation}
    where $
    \Theta_i^{(d)} := X_i^{(d)}/\| X_i^{(d)} \|, \ i=1,2,$
  are the angular components of $X_1^{(d)}$ and
  $X_2^{(d)}$.
\item With the same sequence $(a(n))_{n \in \N}$ as in (I),
  \begin{equation}
    \label{equ_ass_iii}
    \lim_{s \to 0+} \limsup_{d \to \infty} \frac{n}{\sqrt{a(n)}}
    \big\|\E X_1^{(d)} \one_{\| X_1^{(d)} \|^2 \leq s a(n)}\big\|=0.
  \end{equation}
\end{enumerate}

Let
\begin{equation}
  \label{equ:poi_def}
  \sP:=\sum_{k} \delta_{(x_k, y_k)}
\end{equation}
be the Poisson process on 
$[0, \infty) \times (0, \infty) $ with intensity measure 
$\mathbb{LEB} \times \nu$, where $\mathbb{LEB}$ denotes the 
Lebesgue measure and $\nu$ is the L\'evy measure as in (I). Thus
$\zeta_t:=\sum_{x_k\leq t}y_k$ is a subordinator.

Consider a $d$-dimensional random walk with non-integrable
increments satisfying conditions (I)-(III). \cite{Kabluchkoetal2024}
proved that, as $n,d \to \infty$, the random metric
space $\big( a(n)^{-1/2}\sM_{n}^{(d)}, \| \cdot \|\big)$
converges in distribution to the random metric space $\big([0,1],\rho\big)$, where
\begin{displaymath}
  \rho(s,t):=\big|\zeta_t-\zeta_s\big|^{1/2},
  \quad 0\leq s,t\leq 1.
\end{displaymath}

Following a similar approach in \cite{Kabluchkoetal2024},
we can prove the bridge variant of this result. 

\begin{theorem}
  \label{thm_br_cs}
  Under the assumptions (I)-(III), as $n, d \to\infty$, the random
  metric space $\big(a(n)^{-1/2}\sB_{n}^{(d)}, \| \cdot
  \|\big)$ converges in distribution to $\big([0,1],\rho_{B}\big)$,
  where
  \begin{displaymath}
    \rho_{B}(s,t):=\big|\zeta_t-\zeta_s-(t-s)\zeta_1\big|^{1/2},
    \quad 0\leq s,t\leq 1.
  \end{displaymath}
\end{theorem}

As in Remark \ref{rem:1}, the continuous mapping
  theorem is used to verify the validity of Theorem \ref{thm_br_cs},
  noticing that the map \eqref{eq:add_1} is continuous on 
  c{\`a}dl{\`a}g function.

\section{Proof of Theorem~\ref{the_br_l2}}
The following result covers the major step in the proof of
Theorem~\ref{the_br_l2}.

\begin{lemma}
  \label{lem:con_sig}
  Let $n = n(d)$ be an arbitrary sequence of positive integers such
  that $n(d)\to\infty$, as $d\to\infty$. Under the assumptions (a)–(d),
  \begin{equation}
  \label{equ_br_single_l_2}
  n^{-1}\Big\| S_{\lfloor nt \rfloor}^{(d)} -
   \frac{\lfloor nt \rfloor}{n} S_n^{(d)}\Big \|^2  \overset{p} \to
    (1 - t)t \quad\text{as} \ d \to \infty.
  \end{equation}
\end{lemma}

\begin{pf}
  For every $i \in \N$,
\begin{align*}
  \Big \| S_{i}^{(d)} - \frac{i}{n} S_n^{(d)} \Big\|^2
  &=\bigg(1-\frac{i}{n}\bigg)^2\big\| S_i^{(d)}\big\|^2
  + 2\frac{i}{n}\bigg(\frac{i}{n}-1\bigg)
  \langle S_i^{(d)}, S_n^{(d)}-S_i^{(d)}\rangle
    +\bigg(\frac{i}{n}\bigg)^2\big\| S_n^{(d)}-S_i^{(d)}\big\|^2\\
  &=\bigg(1-\frac{i}{n}\bigg)\big\| S_i^{(d)}\big\|^2
    +\frac{i}{n}\big\| S_n^{(d)}-S_i^{(d)}\big\|^2
    -\frac{i}{n}\bigg(1-\frac{i}{n}\bigg)\big\| S_n^{(d)}\big\|^2.
\end{align*}

It is clear that
\begin{align*}
  \bigg(1-\frac{i}{n}\bigg) \big\| S_i^{(d)}\big\|^2
  &=\bigg(1-\frac{i}{n}\bigg)T_i^{(d)}
  +\bigg(1-\frac{i}{n}\bigg)Q_{1,i}^{(d)},\\
  \frac{i}{n}\big\| S_n^{(d)}-S_i^{(d)}\big\|^2
  &=\frac{i}{n}\big(T_n^{(d)}-T_i^{(d)}\big)
  +\frac{i}{n}Q_{i+1,n}^{(d)},\\ \intertext{and}
  \frac{i}{n}\bigg(1-\frac{i}{n}\bigg) \big\| S_n^{(d)}\big\|^2
  &=\frac{i}{n}\bigg(1-\frac{i}{n}\bigg)T_n^{(d)}
  +\frac{i}{n}\bigg(1-\frac{i}{n}\bigg)Q_{1,n}^{(d)},
\end{align*}
where 
\begin{equation}
  \label{equ:TQ}
  T_i^{(d)}:=\sum_{l=1}^{i}\big\|X_l^{(d)}\big\|^2\quad\text{and}\quad
  Q_{a,b}^{(d)}:=\sum_{a\leq l\neq k\leq b}\langle X_l^{(d)},X_k^{(d)}\rangle.
\end{equation}

The weak law of large numbers implies that for all $t\in[0,1]$, as
$d\to\infty$,
\begin{equation}
  \label{equ:single_T}
    n^{-1}\bigg(1-\frac{\nt}{n}\bigg)T_{\nt}^{(d)}
   \overset{p}{\to}(1-t)t,\quad
    n^{-1}\frac{\nt}{n}\big(T_n^{(d)}-T_{\nt}^{(d)}\big)\overset{p}{\to}t(1-t),\quad
    n^{-1}\frac{\nt}{n}\bigg(1-\frac{\nt}{n}\bigg)
    T_n^{(d)}\overset{p}{\to}t(1-t).
\end{equation}

By conditions (a)-(d), for all $a\leq i\neq j\leq b$,
$1\leq a\leq b\leq n$,
\begin{align}
  \label{equ:Q_main}
    \E \big(Q_{a,b}^{(d)}\big)^2
    &=\E \bigg(\sum_{a\leq k\neq l\leq b}
    \sum_{j= 1}^{d} X_{k,j}^{(d)} X_{l,j}^{(d)} \bigg )^2  
    = \E \sum_{a\leq k\neq l\leq b} 
    \sum_{a\leq k'\neq l'\leq b}
    \sum_{j=1}^{d} \sum_{j'=1}^{d} X_{k,j}^{(d)} X_{l,j}^{(d)}
    X_{k',j'}^{(d)} X_{l',j'}^{(d)} \notag\\
    &= 2 \sum_{a\leq k\neq l\leq b} \sum_{j= 1}^{d}
    \E \big( X_{k,j}^{(d)} \big)^2  \E \big( X_{l,j}^{(d)} \big)^2 
    \leq 2(b-a+1)(b-a)\max_{1\leq j\leq d}
    \E (X_{1,j})^2 \sum_{j=1}^{d} \E (X_{1,j})^2\\
    &= 2(b-a+1)(b-a)\max_{1\leq j\leq d}
    \E (X_{1,j})^2\notag.
\end{align} 
Then Markov's inequality implies that
\begin{align*}
  \Prob{Q_{1,\nt}^{(d)}\geq n \eps}
  \leq n^{-2}\eps^{-2}\E (Q_{1,\nt}^{(d)})^2
  &\leq 2 n^{-2}\eps^{-2} (\nt)(\nt-1)\max_{1\leq j\leq d}
  \E (X_{1,j})^2\\
  &\leq 2\eps^{-2} \max_{1\leq j\leq d}
  \E (X_{1,j})^2 \to 0\quad\text{as}\ d\to\infty.
\end{align*}

By the same argument, $Q_{\nt+1,n}^{(d)}\overset{p}{\to}0$ 
and $Q_{1,n}^{(d)}\overset{p}{\to}0$. Finally, 
\eqref{equ_br_single_l_2} is completed by combining
\eqref{equ:single_T}, \eqref{equ:Q_main} and Markov's inequality.
\end{pf}

\begin{pf}[of Theorem \ref{the_br_l2}]
Fix an $m\in\N$. By Lemma \ref{lem:con_sig}, we obtain that, for all
$i=0,\ldots,m$,
\begin{equation}
  \label{equ_point_l2}
  n^{-1/2}\Big\| S_{\nt[i/m]}^{(d)} - \frac{\nt[i/m]}{n} S_n^{(d)}\Big\|
  \overset{p} \to \sqrt{(1 - i/m)i/m}
  \quad \text{as}\ d\to\infty.
\end{equation}
By applying the same method used to prove Lemma \ref{lem:con_sig},
\begin{multline*}
  \max_{0 \leq i \leq j \leq m}\Bigg|n^{-1/2}
  \bigg\| \Big(S_{\nt[j/m]}^{(d)} - \frac{\nt[j/m]}{n} S_n^{(d)}\Big)
    - \Big(S_{\nt[i/m]}^{(d)} - \frac{\nt[i/m]}{n} S_n^{(d)} \Big)\bigg\|\\
  - \sqrt{\big(1 - ( j/m - i/m) \big)}
    \sqrt{j/m - i/m} \Bigg| \overset{p} \to 0
    \quad \text{as}\ d\to\infty.
\end{multline*}

For $s \in [i/m, (i+1)/m]$ and
$t \in [j/m, (j+1)/m]$ with $0 \leq s \leq t \leq 1$,
by the triangle inequality,
\begin{multline*}
  \Bigg| n^{-1/2}\bigg\| \Big(S_{\nt}^{(d)} - \frac{\nt}{n} S_n^{(d)}\Big)
    - \Big( S_{\ns}^{(d)} - \frac{\ns}{n} S_n^{(d)} \Big)\bigg \| \\
  \qquad \qquad \qquad \qquad \qquad
  - n^{-1/2}\bigg\| \Big(S_{\nt[j/m]}^{(d)} - \frac{\nt[j/m]}{n} S_n^{(d)}\Big)
    - \Big(S_{\nt[i/m]}^{(d)} - \frac{\nt[i/m]}{n} S_n^{(d)}\Big)\bigg\|\Bigg| \\
  \leq \sup_{z \in [i/m, (i+1)/m]}
  n^{-1/2}\bigg\| \Big(S_{\nz}^{(d)} - \frac{\nz}{n} S_n^{(d)}\Big)
    - \Big(S_{\nt[i/m]}^{(d)} - \frac{\nt[i/m]}{n} S_n^{(d)}\Big)\bigg\| \\
  + \sup_{z \in [j/m, (j+1)/m]}
  n^{-1/2}\bigg\| \Big(S_{\nz}^{(d)} - \frac{\nz}{n} S_n^{(d)}\Big)
    - \Big(S_{\nt[j/m]}^{(d)} - \frac{\nt[j/m]}{n} S_n^{(d)} \Big)\bigg\|.
\end{multline*}

It suffices to verify that, as $d\to\infty$,
\begin{equation*}
  \P\bigg\{\max_{0\leq i\leq m-1} \sup_{z \in{[i/m, (i+1)/m]}} 
  \bigg\| \Big(S_{\nz}^{(d)} - \frac{\nz}{n} S_n^{(d)}\Big)
    - \Big(S_{\nt[i/m]}^{(d)} - \frac{\nt[i/m]}{n} S_n^{(d)}\Big)\bigg\|^2
  \geq n\eps\bigg\}\to 0.
\end{equation*} 

By the union bound, it suffices to prove that 
\begin{equation}
\label{equ:diff_expand}
  \P\bigg\{\sup_{z\in [i/m, (i+1)/m]}
  \bigg\| \Big(S_{\nz}^{(d)} - \frac{\nz}{n} S_n^{(d)}\Big)
  -\Big(S_{\nt[i/m]}^{(d)} - \frac{\nt[i/m]}{n}
  S_n^{(d)}\Big)\bigg\|^2 \geq n\eps\bigg\}\to 0
\end{equation}
as $d \to \infty$ for all $i=0,\ldots,m-1$.

For all $z\in\big[i/m,(i+1)/m\big]$,
\begin{align*}
  \bigg\| \Big(S_{\nz}^{(d)} &- \frac{\nz}{n} S_n^{(d)}\Big)
   -\Big(S_{\nt[i/m]}^{(d)} - \frac{\nt[i/m]}{n}
   S_n^{(d)}\Big)\bigg\|^2
   =\bigg(1-\frac{\nz}{n}+\frac{\nt[i/m]}{n}\bigg)\big\|S_{\nz}^{(d)}\big\|^2\\
   &+\big\|S_{\nt[i/m]}^{(d)}\big\|^2
   -\Bigg(\bigg(\frac{\nz}{n}-\frac{\nt[i/m]}{n}\bigg)
   -\bigg(\frac{\nz}{n}-\frac{\nt[i/m]}{n}\bigg)^2\Bigg)\big\|S_n^{(d)}\big\|^2\\
   &-2\langle S_{\nz}^{(d)},S_{\nt[i/m]}^{(d)} \rangle
   +\bigg(\frac{\nz}{n}-\frac{\nt[i/m]}{n}\bigg) 
   \big\|S_n^{(d)}-S_{\nz}^{(d)}\big\|^2
   +2\bigg(\frac{\nz}{n}-\frac{\nt[i/m]}{n}\bigg) 
   \langle S_{\nt[i/m]}^{(d)},S_n^{(d)} \rangle.
\end{align*}

For proving \eqref{equ:diff_expand}, it suffices to show that, as $d\to\infty$,
\begin{equation}
  \label{equ:first_part}
    \sup_{z\in[i/m,(i+1)/m]}\Bigg|n^{-1}
    \bigg(1-\frac{\nz}{n}+\frac{\nt[i/m]}{n}\bigg)
    \big\|S_{\nz}^{(d)}\big\|^2-(1-z+i/m)z\Bigg|\overset{p}\to 0,
\end{equation}
\begin{align}
  \label{equ:check_2}
     \sup_{z\in[i/m,(i+1)/m]}\Bigg|n^{-1}\big\|S_{\nt[i/m]}^{(d)}\big\|^2
   &-n^{-1}\Bigg(\bigg(\frac{\nz}{n}-\frac{\nt[i/m]}{n}\bigg)
   -\bigg(\frac{\nz}{n}-\frac{\nt[i/m]}{n}\bigg)^2\Bigg)\big\|S_n^{(d)}\big\|^2\notag\\
   &+n^{-1}2\bigg(\frac{\nz}{n}-\frac{\nt[i/m]}{n}\bigg) 
   \langle S_{\nt[i/m]}^{(d)},S_n^{(d)} \rangle\\
   &-\Big(i/m-\big((z-i/m)-(z-i/m)^2\big)
   +2(z-i/m)i/m\Big)\Bigg|\overset{p}\to 0,\notag
 \end{align}
\begin{equation}
  \label{equ:check_3}
  \sup_{z\in[i/m,(i+1)/m]}
  \Bigg|n^{-1}\bigg(\frac{\nz}{n}-\frac{\nt[i/m]}{n}\bigg) 
  \Big\|S_n^{(d)}-S_{\nz}^{(d)}\Big\|^2
  -(z-i/m)(1-z)\Bigg|\overset{p}\to 0,
\end{equation}
and
\begin{equation}
  \label{equ:last_part}
  \sup_{z\in[i/m,(i+1)/m]}
  \bigg|n^{-1}2\langle S_{\nz}^{(d)},S_{\nt[i/m]}^{(d)} \rangle-2i/m\bigg|
  \overset{p}{\to}0.
\end{equation}

Combining \eqref{equ:first_part}-\eqref{equ:last_part},
the triangle inequality yields that for all $0\leq i\leq m-1$,
\begin{multline*}
  \sup_{z \in{[i/m, (i+1)/m]}} \bigg|n^{-1}
  \Big\| \Big(S_{\nz}^{(d)} - \frac{\nz}{n} S_n^{(d)}\Big)
    - \Big(S_{\nt[i/m]}^{(d)} - \frac{\nt[i/m]}{n} S_n^{(d)}\Big)\Big\|^2\\
  -\big(1-(z-i/m)\big)(z-i/m)
  \bigg|\overset{p}{\to}0\quad\text{as}\ d\to\infty.
\end{multline*}
Therefore, \eqref{equ:diff_expand} holds for every $m\geq 1/\eps$, since
\begin{displaymath}
  \max_{0\leq i\leq m-1}\sup_{z\in[i/m,(i+1)/m]}\big(1-(z-i/m)\big)
  (z-i/m)\leq1/m\leq \eps.
\end{displaymath}

By the triangle inequality, \eqref{equ:first_part} is bounded above by
\begin{displaymath}
  \sup_{z\in[i/m,(i+1)/m]}\Bigg|
    n^{-1}\bigg(1+\frac{\nt[i/m]}{n}\bigg)\big\|S_{\nz}^{(d)}\big\|^2-(1+i/m)z\Bigg|
    + \sup_{z\in[i/m,(i+1)/m]}\Bigg|n^{-1}
    \bigg(\frac{\nz}{n}\bigg)\big\|S_{\nz}^{(d)}\big\|^2-z^2\Bigg|.
\end{displaymath}
Recall notation from \eqref{equ:TQ}.
Then it suffices to show that
\begin{displaymath}
  n^{-1}\bigg(1+\frac{\nt[i/m]}{n}\bigg)\big\|S_{\nz}^{(d)}\big\|^2
  =n^{-1}\bigg(1+\frac{\nt[i/m]}{n}\bigg)T_{\nz}^{(d)}
  +n^{-1}\bigg(1+\frac{\nt[i/m]}{n}\bigg)Q_{1,\nz}^{(d)}\overset{p}{\to}(1+i/m)z,
\end{displaymath}
and
\begin{displaymath}
  n^{-1} \frac{\nz}{n}\big\|S_{\nz}^{(d)}\big\|^2
  =n^{-1}\frac{\nz}{n}T_{\nz}^{(d)}
  +n^{-1}\frac{\nz}{n}Q_{1,\nz}^{(d)}\overset{p}{\to}z^2
\end{displaymath}
uniformly over $z\in\big[i/m,(i+1)/m\big]$. The weak law of large
numbers and Dini's second theorem (see Exercise 127 on page 81 and its
proof on page 270 from \cite{PolyaSzego1998}) imply that
\begin{displaymath}
  \sup_{z\in[i/m,(i+1)/m]}\Bigg|n^{-1}
  \bigg(1+\frac{\nt[i/m]}{n}\bigg)T_{\nz}^{(d)}-(1+i/m)z\Bigg|
  \overset{p}{\to}0\quad\text{as}\ d\to\infty,
\end{displaymath}
and 
\begin{displaymath}
  \sup_{z\in[i/m,(i+1)/m]}\Bigg|n^{-1}
  \frac{\nz}{n}T_{\nz}^{(d)}-z^2\Bigg|
  \overset{p}{\to}0\quad\text{as}\ d\to\infty.
\end{displaymath}

Note that
\begin{displaymath}
  \P\Bigg\{\sup_{z\in[i/m,(i+1)/m]}\Bigg|
  \bigg(1+\frac{\nt[i/m]}{n}\bigg)Q_{1,\nz}^{(d)} \Bigg| \geq n \eps\Bigg\}
  \leq \P\bigg\{\sup_{z\in[i/m,(i+1)/m]} 
  \big|Q_{1,\nz}^{(d)}\big| \geq n \eps/2\bigg\},
\end{displaymath}
and
\begin{displaymath}
  \P\Bigg\{\sup_{z\in[i/m,(i+1)/m]}\Bigg|
    \frac{\nz}{n} Q_{1,\nz}^{(d)}\Bigg| \geq n \eps\Bigg\}
  \leq \P\bigg\{\sup_{z\in[i/m,(i+1)/m]}
   \big|Q_{1,\nz}^{(d)}\big| \geq n \eps/2\bigg\}.
\end{displaymath}
The sequence
$$Q_{1,i}^{(d)}= 2 \sum_{k=1}^{i} \langle X_k^{(d)}, S_{k - 1}^{(d)}
\rangle,\quad i\in\N,$$ is a martingale with respect to the filtration
$\big\{\salg_{i}^{(d)} = \sigma(X_1^{(d)}, \ldots, X_i^{(d)}), i \geq 1
\big\}$, since
\begin{displaymath}
  \E \big(Q_{1,i}^{(d)} - Q_{1,i-1}^{(d)} \mid \salg_{i - 1}^{(d)}\big)
  = 2 \E \big(\langle X_i^{(d)}, S_{i - 1}^{(d)} \rangle \mid\salg_{i - 1}^{(d)}\big)
  = 2 \sum_{j = 1}^{d} \E \big(X_{i,j}^{(d)} S_{i - 1,j}^{(d)} 
  \mid  \salg_{i - 1}^{(d)}\big) \overset{\text{a.s.}}{=} 0.
\end{displaymath}
By Doob's inequality, conditions (a)-(d) and \eqref{equ:Q_main},
\begin{multline*}
  \P\bigg\{\sup_{z\in[i/m,(i+1)/m]} \big | Q_{1,\nz}^{(d)}\big |
  \geq n \eps/2\bigg\}
  \leq n^{-2}\eps^{-2}\E\big(Q_{1,\nt[(i+1)/m]}^{(d)}\big)^2\\
  \leq 8n^{-2}\eps^{-2}\nt[(i+1)/m](\nt[(i+1)/m]-1)
  \max_{1\leq j\leq d}\E(X_{1,j}^{(d)})\to 0.
\end{multline*} 

The term given by \eqref{equ:check_2} is bounded above by
\begin{align*}
    n^{-1}&\bigg|\big\|S_{\nt[i/m]}^{(d)}\big\|^2-i/m\bigg|
   +\sup_{z\in[i/m,(i+1)/m]}
   \Bigg|n^{-1}2\Bigg(\frac{\nz}{n}-\frac{\nt[i/m]}{n}\bigg) 
   \langle S_{\nt[i/m]}^{(d)},S_n^{(d)}
   \rangle-2\Big(z-i/m\Big)i/m\Bigg|\\
   &\qquad\qquad\qquad\qquad\qquad\qquad
   +\sup_{z\in[i/m,(i+1)/m]}
   \Bigg|n^{-1}\bigg(\frac{\nz}{n}-\frac{\nt[i/m]}{n}\bigg)\big\|S_n^{(d)}\big\|^2
   -\big(z-i/m\big)\Bigg|\\
   &\qquad\qquad\qquad\qquad\qquad\qquad
   +\sup_{z\in[i/m,(i+1)/m]}
   \Bigg|n^{-1}\bigg(\frac{\nz}{n}-\frac{\nt[i/m]}{n}\bigg)^2\big\|S_n^{(d)}\big\|^2
   -\big(z-i/m\big)^2\Bigg|.
\end{align*}
It is clear that
\begin{displaymath}
  \big\|S_{\nt[i/m]}^{(d)}\big\|^2=T_{\nt[i/m]}^{(d)}+Q_{1,\nt[i/m]}^{(d)},\quad
   \langle S_{\nt[i/m]}^{(d)},S_n^{(d)} \rangle
  =T_{\nt[i/m]}^{(d)}
  +\sum_{l=1}^{\nt[i/m]}\sum_{1\leq k\leq n,k\neq l}
  \langle X_l^{(d)},X_k^{(d)}\rangle,\quad
  \big\|S_n^{(d)}\big\|^2=T_n^{(d)}+Q_{1,n}^{(d)}.
\end{displaymath}
The weak law of large numbers and Dini's second theorem imply
that, as $d\to\infty$,
\begin{align*}
  n^{-1}T_{\nt[i/m]}^{(d)}&\overset{p}{\to} i/m,\quad
  n^{-1}2\bigg(\frac{\nz}{n}-\frac{\nt[i/m]}{n}\bigg) 
  T_{\nt[i/m]}^{(d)}
  \overset{p}{\to} 2\Big(z-i/m\Big)i/m,\\\intertext{and}
  n^{-1}\bigg(\frac{\nz}{n}-\frac{\nt[i/m]}{n}\bigg)T_{n}^{(d)}
  &\overset{p}{\to}\big(z-i/m\big),\quad
  n^{-1}\bigg(\frac{\nz}{n}-\frac{\nt[i/m]}{n}\bigg)^2T_{n}^{(d)}
  \overset{p}{\to}\big(z-i/m\big)^2,
\end{align*}
uniformly over $z\in[i/m,(i+1)/m]$. By Markov's
inequality and \eqref{equ:Q_main},
\begin{displaymath}
  \Prob{Q_{1,\nt[i/m]}^{(d)}\geq n\eps}
  \leq 2n^{-2}\eps^{-2}\nt[i/m]\big(\nt[i/m]-1\big) 
  \max_{1\leq j\leq d}\E(X_{1,j}^{(d)})\to 0,
\end{displaymath}
\begin{multline*}
  \P\bigg\{\sup_{z\in[i/m,(i+1)/m]}
    2\bigg(\frac{\nz}{n}-\frac{\nt[i/m]}{n}\bigg)
    \Big|\sum_{l=1}^{\nt[i/m]}\sum_{1\leq k\leq n,k\neq l}
    \langle X_l^{(d)},X_k^{(d)}\rangle\Big|\geq n\eps\bigg\}\\
  \leq 4n^{-2}\eps^{-2}\E\bigg(\sum_{l=1}^{\nt[i/m]}\sum_{1\leq k\leq n,k\neq l}
    \langle X_l^{(d)},X_k^{(d)}\rangle\bigg)^2\to 0,
\end{multline*}
\begin{displaymath}
  \P\bigg\{\sup_{z\in[i/m,(i+1)/m]}
    \bigg(\frac{\nz}{n}-\frac{\nt[i/m]}{n}\bigg)
    \Big|Q_{1,n}^{(d)}\Big|\geq n\eps\bigg\}
  \leq 2n^{-2}\eps^{-2}n(n-1)
  \max_{1\leq j\leq d}\E(X_{1,j}^{(d)})\to 0,
\end{displaymath}
and
\begin{displaymath}
  \P\bigg\{\sup_{z\in[i/m,(i+1)/m]}
   \bigg(\frac{\nz}{n}-\frac{\nt[i/m]}{n}\bigg)^2
   \Big|Q_{1,n}^{(d)}\Big|\geq n\eps\bigg\}
  \leq 2n^{-2}\eps^{-2}n(n-1)
   \max_{1\leq j\leq d}\E(X_{1,j}^{(d)})\to 0,
\end{displaymath}
which shows that \eqref{equ:check_2} converges uniformly to 0 in
probability as $d\to\infty$.

By the triangle inequality, \eqref{equ:check_3}
is bounded above by
\begin{multline*}
  \sup_{z\in[i/m,(i+1)/m]}
  \bigg|n^{-1}\bigg(\frac{\nz}{n}-\frac{\nt[i/m]}{n}\bigg) 
  \big(T_n^{(d)}-T_{\nz}^{(d)}\big)
   -\Big(z-i/m\Big)(1-z)\bigg|\\
   +\sup_{z\in[i/m,(i+1)/m]}
   \bigg|n^{-1}\sum_{\nz+1\leq l\neq k\leq n} 
   \langle X_l^{(d)},X_k^{(d)} \rangle\bigg|\\
   \leq \sup_{z\in[i/m,(i+1)/m]}
  \bigg|n^{-1}\frac{\nz}{n}T_n^{(d)}-z\bigg|
   +\sup_{z\in[i/m,(i+1)/m]}
   \bigg|n^{-1}\frac{\nz}{n}T_{\nz}^{(d)}-z^2\bigg|\\
   +\sup_{z\in[i/m,(i+1)/m]}
  \bigg|n^{-1}\frac{\nt[i/m]}{n}
  \big(T_n^{(d)}-T_{\nz}^{(d)}\big)
   -i/m(1-z)\bigg|
   + \sup_{z\in[i/m,(i+1)/m]}
   \bigg|n^{-1}Q_{\nz+1,n}^{(d)}\bigg|.
\end{multline*}
The first three terms converge to 0 by the weak law of large numbers
and Dini's second theorem. It is clear that for all
$i \in \N$, $$Q_{i+1,n}^{(d)}
= 2 \sum_{k=i+1}^{n} \langle X_{k}^{(d)},S_n^{(d)}-S_{k}^{(d)} \rangle$$
is a reverse martingale with respect to the filtration
$\big\{\salg_{i+1}'^{(d)} = \sigma(X_{i+1}^{(d)}, \ldots, X_n^{(d)}), i \geq 1 \big\}$, since
\begin{displaymath}
  \E \Big(Q_{i+1,n}^{(d)} - Q_{i+2,n}^{(d)} \mathrel{\big|} \salg_{i +2}'^{(d)}\Big)
  = 2 \E \Big(\langle X_{i+1}^{(d)}, S_{n}^{(d)}-S_{i+1}^{(d)} \rangle
  \mathrel{\big|}\salg_{i+2}'^{(d)}\Big)\notag
  = 2 \sum_{j = 1}^{d} \E \Big(X_{i+1,j}^{(d)}
  \big(S_{n,i}^{(d)}-S_{i+1,j}^{(d)}\big)
  \mathrel{\big|}  \salg_{i+2}'^{(d)}\Big) \overset{\text{a.s.}}{=}0.
\end{displaymath}
Doob's inequality and \eqref{equ:Q_main} yield that for all $1\leq i\leq m-1$,
\begin{align*}
  \P\bigg\{\sup_{z\in[i/m,(i+1)/m]}\Big|
  Q_{\nz+1,n}^{(d)} \Big|\geq n\eps\bigg\}
  &\leq n^{-2}\eps^{-2}\E Q_{\nt[i/m]+1,n}^{(d)})^2\\
  &\leq 2n^{-2}\eps^{-2}(n-\nt[i/m])(n-\nt[i/m]-1)\max_{1\leq i\leq d}
    \E (X_{1,i})^2\to 0.
\end{align*}

It is left to show that \eqref{equ:last_part} holds. By the triangle
inequality, this term is bounded above by
\begin{equation}
  \label{equ:last_term_bound}
  \sup_{z\in[i/m,(i+1)/m]}
  \bigg|2T_{\nt[i/m]}^{(d)}-2i/m\bigg|
  +\sup_{z\in[i/m,(i+1)/m]}
  \bigg|2\sum_{k=1}^{\nt[i/m]}\sum_{1\leq l\leq\nz, k\neq l}
  \langle X_k^{(d)},X_l^{(d)}\rangle\bigg|.
\end{equation}
The first term converges to 0 by the weak law of large numbers and
Dini's second theorem. Note that
\begin{displaymath}
\sum_{k=1}^{\nt[i/m]}\sum_{1\leq l\leq r, l\neq k} \langle
X_k^{(d)},X_l^{(d)}\rangle,\quad r=\nt[i/m],\ldots,n,
\end{displaymath}
is a martingale with respect to the
filtration $\big\{\salg_{r}^{(d)} = \sigma(X_1^{(d)},\ldots, X_{r}^{(d)}), r
\geq 1 \big\}$, since
\begin{align*}
  \E\bigg(\sum_{k=1}^{\nt[i/m]}\sum_{1\leq l\leq r+1, l\neq k} &\langle
  X_k^{(d)},X_l^{(d)}\rangle
  -\sum_{k=1}^{\nt[i/m]}\sum_{1\leq l\leq r, l\neq k} \langle
  X_k^{(d)},X_l^{(d)}\rangle\mathrel{\big|}\salg_{r} \bigg)\\
  &=\E\bigg(\sum_{k=1}^{\nt[i/m]}\langle
  X_k^{(d)},X_{r+1}^{(d)}\rangle\mathrel{\big|}\salg_{r}\bigg)
  =\sum_{j=1}^{d}\sum_{k=1}^{\nt[i/m]}X_{k,j}^{(d)}
  \E X_{r+1,j}^{(d)}=0.
\end{align*}
The second term of \eqref{equ:last_term_bound} converges to 0 
by Doob's inequality, namely,
\begin{align*}
  \P\bigg\{\sup_{z\in[i/m,(i+1)/m]}
  &  \Big|\sum_{k=1}^{\nt[i/m]}\sum_{1\leq l\leq \nz, l\neq k}
  \langle X_k^{(d)},X_l^{(d)}\rangle \Big|\geq n\eps/2\bigg\}\\
  &\leq 4n^{-2}\eps^{-2} \E \bigg(
    \sum_{k=1}^{\nt[i/m]}\sum_{1\leq l\leq \nt[(i+1)/m],l\neq k}
    \sum_{i= 1}^{d} X_{k,i}^{(d)} X_{l,i}^{(d)} \bigg )^2  \\
  &= 8n^{-2}\eps^{-2}\sum_{k=1}^{\nt[i/m]}\sum_{1\leq l\leq \nt[(i+1)/m],l\neq k}
    \sum_{i=1}^{d}\E \big( X_{k,i}^{(d)} \big)^2  \E \big( X_{l,i}^{(d)} \big)^2 \\
  & \leq 8n^{-2}\eps^{-2}\nt[i/m](\nt[(i+1)/m]- 1)\max_{1\leq i\leq d}
    \E (X_{1,i})^2 \sum_{i = 1}^{d} \E (X_{1,i})^2 \to 0,
\end{align*}
which completes the proof.
\end{pf}

\section{Proof of Theorem \ref{thm_br_cs}}

Instead of proving the weak convergence to $\big([0,1],\rho_{B}\big)$,
we substitute the limit with an isometric space constructed as
follows. Recall the Poisson process $\mathcal{P}$ from
\eqref{equ:poi_def} and define 
\begin{displaymath}
  \sR_{\nu} :=\mathrm{cl} \bigg\{ \sum_{k : x_k \leq t} e_k y_k^{1/2}
  -t\sum_{k: x_k \leq 1} e_k y_k^{1/2}: 0\leq t\leq 1\bigg\},
\end{displaymath}
where $(e_k)_{k\in\N}$ is the standard orthonormal basis of the
space $\ell^2$ of square-summable sequences. Then $\sR_{\nu}$
considered as a compact subset of $\ell^2$ is isometric to
$\big([0,1],\rho_{B}\big)$. 

Next, we construct a further metric space which approximates
$(\sR_{\nu},\|\cdot\|)$.  For $s > 0$, define the truncated
variables
\begin{equation*}
  X_k^{(d)} (s) := X_k^{(d)} \one_{\|X_k^{(d)} \|^2 \geq s a(n)},
  \quad k \in \N,
\end{equation*}
the corresponding random walks
\begin{equation*}
    S_0^{(d)}(s) := 0, \qquad S_k^{(d)} (s) := X_1^{(d)} (s)+\cdots
    +X_k^{(d)} (s),
\end{equation*}
and the sets
\begin{equation*}
  \sB_k^{(d)} (s) :=\Big\{S_0^{(d)}(s),S_1^{(d)}(s)-\frac{1}{k} S_k^{(d)}(s),
  \ldots, S_k^{(d)}(s)-\frac{k}{k} S_k^{(d)}(s)\Big\}.
\end{equation*}
Also define the random set
\begin{equation*}
  \sR_{\nu} (s) := \Big\{ \sum_{k: x_k \leq t} e_ky_k^{1/2}\one_{y_k > s}
  -t\sum_{k: x_k \leq 1} e_k y_k^{1/2}
  \one_{y_k > s}, 0 \leq t \leq 1 \Big\}.
\end{equation*}

It suffices to check that the following three statements hold true.
\begin{enumerate}[(i)]
\item For every fixed $s > 0$, in the Gromov-Hausdorff sense,
  $\big(a_n^{-1/2} \sB_{n}^{(d)}(s), \| \cdot \| \big) \overset{d}
  \to \big( \sR_{\nu} (s), \| \cdot \| \big)\,\text{as} \ d \to\infty.$
\item Under the Gromov-Hausdorff distance,
  $\big(  \sR_{\nu} (s), \| \cdot \| \big) \overset{a.s.}{\to}
  \big(\sR_{\nu}, \| \cdot \| \big) \,\text{as} \ s \downarrow 0.$
\item For every fixed $\eps > 0$,
  $\lim_{s \to 0+} \limsup_{d \to \infty}
  \Prob{\dGH \Big(a_n^{-1/2}\sB_{n}^{(d)}(s), a_n^{-1/2}\sB_{n}^{(d)} \Big) > \eps}=0.$
\end{enumerate}

First, we prove (i).
For $k \in \N$, set $R_k^{(d)} := \| X_k^{(d)} \|^2$. Under 
\eqref{equ_ass_i_1}, 
\begin{equation*}
  \sP_n:=\sum_{k \geq 0} \delta_{\big(k / n, R_k^{(d)} /a(n)\big)}\overset{d}{\to}
  \sP=\sum_{k\geq 0} \delta_{(x_k, y_k)} \quad\text{as} \ d \to \infty.
\end{equation*}
By Skorokhod's representation theorem, we are allowed to pass the distributional
convergence to the a.s. convergence on a new probability space
$(\Bar{\Omega}, \Bar{\salg}, \mathbf{\Bar{P}})$, which contains
a distributional copy 
$(\Bar{X}_k^{(d)})_{k \in \N}$ of the sequence $(X_k^{(d)})_{k\in\N}$
for every $d \in \N$ and a distributional copy 
$\Bar{\sP} :=\sum_{k\geq 0}\delta_{(\Bar{x}_k, \Bar{y}_k)}$ of the Poisson point process $\sP$,
such that
\begin{equation*}
  \mathcal{\Bar{P}}_n := \sum_{k \geq 0}
    \delta_{\big(k/n,\Bar{R}_k^{(d)}/a(n)\big)} \to
    \mathcal{\Bar{P}}=\sum_{k\geq 0} \delta_{(\Bar{x}_k, \Bar{y}_k)}  \quad
    \mathbf{\Bar{P}}-\text{a.s.} \quad\text{as} \ d \to \infty,
\end{equation*}
where $\Bar{R}_k^{(d)}=\| \Bar{X}_k^{(d)} \|^2$.
Define
\begin{equation*}
  {\Bar{\sB}}_k^{(d)} (s) :=\bigg\{\sum_{j=1}^{l}
  \Bar{X}_j^{(d)}\one_{\Bar{R}_j^{(d)} \geq sa(n)}-\frac{l}{k}
  \sum_{j=1}^{k}\Bar{X}_j^{(d)}\one_{\Bar{R}_j^{(d)} \geq sa(n)},
  l=0, \ldots, k \bigg\}, \quad k \in \N,
\end{equation*}
\begin{equation*}
  \mathcal{\Bar{R}}_{\nu}(s) :=
  \bigg\{ \sum_{k: \Bar{x}_k \leq t}e_k \Bar{y}_k^{1/2}
  \one_{\Bar{y}_k > s}
  -t \sum_{k: \Bar{x}_k \leq 1}
  e_k \Bar{y}_k^{1/2}\one_{\Bar{y}_k > s}, 0 \leq t \leq 1\bigg\}.
\end{equation*}
Since convergence in probability implies convergence in distribution,
it suffices to show that
\begin{equation*}
  \dGH \big(a_n^{-1/2}{\Bar{\sB}}_{n}^{(d)} (s),
  \mathcal{\Bar{R}}_{\nu}(s)\big)
  \overset{\mathbf{\Bar{P}}} \to 0 \quad\text{as} \ d \to \infty.
\end{equation*}
By \citep[Proposition 3.17]{Resnick2007}, for all $\Bar{\omega} \in \Bar{\Omega}$,
there exists an integer $P=P(\Bar{\omega}) \in \N$ and
an enumeration of the atoms of $\mathcal{\Bar{P}}$ and
$\mathcal{\Bar{P}}_n$ in $[0, 1] \times [s, \infty)$ such that for all
sufficiently large $n \in \N$, the restriction of
$\mathcal{\Bar{P}}_n$ onto $\big([0, 1] \times [s, \infty)\big)$ is a Poisson
point process composed by $\big(k_j(n)/n,
\Bar{R}_{k_j(n)}^{(d)}/a(n)\big)$ and the restriction of
$\mathcal{\Bar{P}}$ onto $\big([0, 1] \times [s, \infty)\big)$ is a Poisson
point process composed by $(\Bar{x}_{k_j},\Bar{y}_{k_j})$. Moreover,
\begin{equation*}
    \Big(k_j(n)/n,\Bar{R}_{k_j(n)}^{(d)}/a(n) \Big)
    \overset{\text{a.s.}}{\to}(\Bar{x}_{k_j}, \Bar{y}_{k_j}) \quad\text{as}\ d \to \infty,\
    \text{for all} \ 1\leq j\leq P,\, j\in\N.
\end{equation*}
Without loss of generality, assume that
$\Bar{x}_{k_1}<\Bar{x}_{k_2} < \cdots < \Bar{x}_{k_P}$. Then,
\begin{equation*}
  {\Bar{\sB}}_{n}^{(d)} (s)=\bigg\{
  \sum_{j=1}^{l}\Bar{X}_{k_j(n)}^{(d)}
  -\frac{l}{P} \sum_{j=1}^{P}\Bar{X}_{k_j(n)}^{(d)}, l=0, \ldots, P \bigg\},\quad
  \mathcal{\Bar{R}}_{\nu}(s)=\bigg\{ \sum_{j=1}^{l}
  e_{k_j}\Bar{y}_{k_j}^{1/2}
  -\frac{l}{P} \sum_{j=1}^{P} e_{k_j} \Bar{y}_{k_j}^{1/2},
  l=0,\ldots, P \bigg\}.
\end{equation*}
Note that
\begin{align*}
  \dGH &\big( a_n^{-1/2}\Bar{\sB}_{n}^{(d)} (s), \Bar{\sR}_{\nu}(s)\big) \\
  &\leq 2 \sup_{0 \leq l \leq m \leq P} \Bigg| a_n^{-1/2}
  \bigg\|\sum_{j=l + 1}^{m} \Bar{X}_{k_j(n)}^{(d)}
  -\frac{m-l}{P} \sum_{j=1}^{P}\Bar{X}_{k_j(n)}^{(d)} \bigg\|
  -\bigg\| \sum_{j=l + 1}^{m} e_{k_j} \Bar{y}_{k_j}^{1/2}
  -\frac{m-l}{P} \sum_{j=1}^{P} e_{k_j} \Bar{y}_{k_j}^{1/2}\bigg\| \Bigg|.
\end{align*}
It suffices to verify that, as $d\to\infty$,
\begin{align}
  \label{equ:second_proof_bound_1}
    &\sup_{0 \leq l \leq m \leq P}\Bigg| a_n^{-1/2}\bigg\|
    -\frac{m-l}{P} \bigg( \Bar{X}_{k_1(n)}^{(d)} +
    \cdots+\Bar{X}_{k_l(n)}^{(d)}
    + \Bar{X}_{k_{m +1}(n)}^{(d)} + \cdots
    +\Bar{X}_{k_P(n)}^{(d)}\bigg)
    + \Big(1-\frac{m-l}{P} \Big)\Bar{y}_{k_1}^{1/2}
    \sum_{j=l + 1}^{m} \Bar{X}_{k_j(n)}^{(d)}\bigg\|^2\notag\\
    &-\bigg\| -\frac{m-l}{P} \bigg( e_{k_1} \Bar{y}_{k_1}^{1/2}
    +\cdots + e_{k_l} \Bar{y}_{k_l}^{1/2}
    +e_{k_{m + 1}}\Bar{y}_{k_{m+1}}^{1/2}
    +\cdots +e_{k_P}\Bar{y}_{k_P}^{1/2} \bigg)\\
    &\hspace{8cm}+\Big( 1-\frac{m -l}{P} \Big) \sum_{j=l + 1}^{m}
    e_{k_j} \Bar{y}_{k_j}^{1/2}\bigg\|^2 \Bigg|
    \overset{\mathbf{\Bar{P}}}{\to}0\notag.
\end{align}
The triangle inequality yields that \eqref{equ:second_proof_bound_1}
is bounded above by
\begin{align*}
  \sup_{0 \leq l \leq m \leq P} &\Bigg| \bigg(\frac{m-l}{P}\bigg)^2
      \sum_{j=1}^{l} \bigg(a_n^{-1}\big\|\Bar{X}_{k_j(n)}^{(d)}\big\|^2
      -\Bar{y}_{k_j} \bigg)
      + \bigg( 1-\frac{m-l}{P} \bigg)^2
      \sum_{j=l + 1}^{m} \bigg(a_n^{-1}\big\|\Bar{X}_{k_j(n)}^{(d)}\big\|^2
    -\Bar{y}_{k_j} \bigg) \\
  & + \bigg(\frac{m-l}{P}\bigg)^2
      \sum_{j=m + 1}^{P}\bigg( a_n^{-1}\big\|\Bar{X}_{k_j(n)}^{(d)}\big\|^2
    -\Bar{y}_{k_j} \bigg) \Bigg| \\
  &+ a_n^{-1}\sup_{0 \leq l \leq m \leq P}\Bigg|
      \Big( 1-\frac{m-l}{P} \Big)^2
      \sum_{l+1\leq j\neq j'\leq m } \big\langle\Bar{X}_{k_j(n)}^{(d)},
    \Bar{X}_{k_{j'}(n)}^{(d)} \big\rangle\Bigg| \\
  & + a_n^{-1}\sup_{0 \leq l \leq m \leq P}\Bigg|
      \Big(1-\frac{m-l}{P}\Big) \frac{m-l}{P}
      \sum_{l+1\leq j\leq m, 1\leq j'\leq l} \big\langle
      \Bar{X}_{k_j(n)}^{(d)},\Bar{X}_{k_{j'}(n)}^{(d)}\big\rangle\Bigg|\\
  & + a_n^{-1}\sup_{0 \leq l \leq m \leq P}\Bigg|
      \Big(1-\frac{m-l}{P}\Big) \frac{m-l}{P}
      \sum_{l+1\leq j\leq m, m+1\leq j'\leq P } \big\langle
      \Bar{X}_{k_j(n)}^{(d)},\Bar{X}_{k_{j'}(n)}^{(d)}\big\rangle\Bigg|\\
  & + a_n^{-1}\sup_{0 \leq l \leq m \leq P}\Bigg|
      \Big(\frac{m-l}{P}\Big)^2\sum_{1\leq j\neq j'\leq l }
      \big\langle \Bar{X}_{k_j(n)}^{(d)},
    \Bar{X}_{k_{j'}(n)}^{(d)}\big\rangle\Bigg|\\
  &+ a_n^{-1}\sup_{0 \leq l \leq m \leq P}\Bigg| \Big(\frac{m-l}{P}\Big)^2
      \sum_{m+1\leq j\neq j'\leq P }
      \big\langle \Bar{X}_{k_j(n)}^{(d)},
    \Bar{X}_{k_{j'}(n)}^{(d)}\big\rangle\Bigg|.
\end{align*}
The first term of the right-hand side converges to $0$ for all
$\Bar{\omega} \in \Bar{\Omega}'$. The sum of the last five terms is 
bounded above by
\begin{displaymath}
  5 P^2 a_n^{-1}\sup_{1\leq j\neq j'\leq P} \Big|
  \big\langle\Bar{X}_{k_j(n)}^{(d)},
  \Bar{X}_{k_{j'}(n)}^{(d)} \big\rangle\Big|\leq
  5P^2 a_n^{-1}\sup_{1\leq j\neq j'\leq  n}
  \big| \langle \Bar{X}_j^{(d)} \one_{\|\Bar{X}_j^{(d)} \|^2 \geq s a(n)},
  \Bar{X}_{j'}^{(d)}\one_{\| \Bar{X}_{j'}^{(d)} \|^2 \geq s a(n)}\rangle \big|.
\end{displaymath}

Recall that $\Bar{X}_j^{(d)}$ and $X_j^{(d)}$ have the same
distribution. Hence, it is sufficient to show that
\begin{equation*}
  a_n^{-1}\sup_{1\leq j\neq j'\leq  n}
  \big| \langle X_j^{(d)} \one_{\| X_j^{(d)} \|^2 \geq s a(n)},
  X_{j'}^{(d)} \one_{\| X_{j'}^{(d)}\|^2 \geq sa(n)} \rangle \big|
  \overset{P} \to 0 \quad\text{as} \ d\to\infty.
\end{equation*}
Note that for every fixed $\eps > 0$,
\begin{align*}
  \P&\bigg\{\sup_{1\leq j\neq j'\leq  n}
      \big| \langle X_j^{(d)} \one_{\| X_j^{(d)} \|^2 \geq sa(n)},
      X_{j'}^{(d)} \one_{\| X_{j'}^{(d)} \|^2\geq s a(n)} \rangle\big| \geq \eps a(n)
      \bigg\} \\
  & \leq n^2 \Prob{\big| \langle X_1^{(d)} \one_{\|X_1^{(d)}\|^2 \geq s a(n)},
      X_{2}^{(d)} \one_{\| X_{2}^{(d)} \|^2\geq s a(n)} \rangle \big|\geq \eps a(n)} \\
  & =n^2 \Prob{\big| \langle X_1^{(d)}, X_{2}^{(d)}\rangle \big|>\eps a(n),
      \| X_1^{(d)} \|^2 \geq s a(n), \| X_2^{(d)}\|^2 \geq s a(n),
    \| X_1^{(d)} \| \| X_2^{(d)} \|\geq ra(n)}\\
  & + n^2 \Prob{\big| \langle X_1^{(d)}, X_{2}^{(d)}\rangle\big| >\eps a(n),
      \| X_1^{(d)} \|^2 \geq s a(n), \| X_2^{(d)} \|^2 \geq s a(n),
    \| X_1^{(d)} \| \| X_2^{(d)} \| < ra(n)},
\end{align*}
where $r > s$ is a constant. Note that by \eqref{equ_ass_i_1},
\begin{multline*}
  \lim_{d \to \infty} n^2 \Prob{\| X_1^{(d)} \|^2 \geq s a(n),
    \| X_2^{(d)} \|^2 \geq s a(n), \| X_1^{(d)} \| \| X_2^{(d)}\|\geq r a(n)} \\
  =(\nu \otimes \nu)\big([s, \infty) \times [s, \infty) \cup \{(x,y) \in
  (0, \infty)^2: xy \geq r^2 \}\big).
\end{multline*}
The right-hand side converges to 0 as $r \to +\infty$.
  
Under the assumption \eqref{equ_ass_ii},
\begin{align*}
  n^2 &\Prob{\big| \langle X_1^{(d)}, X_{2}^{(d)}\rangle \big| > \eps a(n),
      \| X_1^{(d)} \|^2 \geq s a(n), \| X_2^{(d)} \|^2\geq s a(n),
    \| X_1^{(d)} \| \| X_2^{(d)} \| < r a(n)}\\
  & \leq n^2 \Prob{\| X_1^{(d)} \|^2 \geq s a(n),
      \|X_2^{(d)} \|^2 \geq s a(n), |
    \langle\Theta_1^{(d)},\Theta_2^{(d)}\rangle | > \eps r^{-1}} \\
  & \leq C(s) \Prob{|\langle \Theta_1^{(d)}, \Theta_2^{(d)}\rangle |>\eps r^{-1}
      \mid \| X_1^{(d)} \|^2 \geq s a(n), \|X_2^{(d)} \|^2 \geq s a(n)} \to 0,
\end{align*}
where $C(s) > 0$ is a constant.

Next we prove (ii). It is known that $\dGH (\sR_{\nu}(s), \sR_{\nu})
\leq\dH(\sR_{\nu}(s), \sR_{\nu})$.
Then it is sufficient to verify that 
$\dH (\sR_{\nu}(s), \sR_{\nu}) \to 0\, \text{as}\ s\downarrow 0.$
\enlargethispage{\baselineskip}
Note that, as $s\downarrow 0$,
\begin{align*}
  \dH (\sR_{\nu}(s), \sR_{\nu})
  & \leq \sup_{t \in [0, 1]} \Big\| \sum_{k: x_k \leq t}e_ky_k^{1/2}\one_{y_k \leq s}
      -t \sum_{k: x_k \leq 1} e_k y_k^{1/2} \one_{y_k \leq s} \Big\| \\
  & \leq \sup_{t \in [0,1]} \big\| \sum_{k: x_k \leq t} e_ky_k^{1/2}\one_{y_k \leq s} \big\|
      + \big\| \sum_{k: x_k \leq 1} e_ky_k^{1/2} \one_{y_k \leq s}\big\| \\
  &=\sup_{t \in [0,1]} \bigg(\sum_{k: x_k \leq t} y_k \one_{y_k \leq s}\bigg)^{1/2}
    + \bigg(\sum_{k: x_k \leq 1} y_k \one_{y_k \leq s}\bigg)^{1/2}
    \leq 2 \bigg( \sum_{k: x_k \leq 1} y_k \one_{y_k \leq s}\bigg)^{1/2}
      \to 0
\end{align*}
by the Lebesgue dominated convergence theorem.

Finally, we prove (iii). By
\citep[Remark 7.3.12]{Buragoetal2001},
\begin{multline*}
  \dGH \Big(a_n^{-1/2}\sB_{ n}^{(d)}(s), a_n^{-1/2}\sB_{ n}^{(d)} \Big)
    \leq \dH\Big(a_n^{-1/2}\sB_{ n}^{(d)}(s), a_n^{-1/2}\sB_{ n}^{(d)}\Big) \\
    \leq a_n^{-1/2} \max_{1\leq k\leq  n}
    \Bigg\|\sum_{j=1}^{k} X_j^{(d)} \one_{\| X_j^{(d)} \|^2 \leq sa(n)}
    -\frac{k}{ n} \sum_{j=1}^{ n} X_j^{(d)} \one_{\|X_j^{(d)} \|^2
      \leq s a(n)} \Bigg\|.
\end{multline*}
Furthermore, the triangle inequality implies that
\begin{displaymath}
  \max_{1\leq k\leq  n} \Bigg\| \sum_{j=1}^{k}X_j^{(d)}\one_{\| X_j^{(d)} \|^2 \leq s a(n)}
      -\frac{k}{ n} \sum_{j=1}^{ n} X_j^{(d)} \one_{\| X_j^{(d)}\|^2 \leq s a(n)}\Bigg\| 
  \leq Z_1^{(d)} (s) + Z_2^{(d)} (s),
\end{displaymath}
where
\begin{align*}
  Z_1^{(d)} (s) &:= \max_{1\leq k\leq  n} \Bigg\| \sum_{j=1}^{k}
       \bigg( X_j^{(d)} \one_{\| X_j^{(d)} \|^2 \leq sa(n)}
    -\E \Big(X_j^{(d)} \one_{\| X_j^{(d)} \|^2 \leq sa(n)}\Big)\bigg) \\
  & \qquad \qquad \qquad \qquad \qquad
       -\frac{k}{ n} \sum_{j=1}^{ n}
       \bigg( X_j^{(d)} \one_{\|X_j^{(d)} \|^2 \leq s a(n)}
       -\E \Big(X_j^{(d)} \one_{\| X_j^{(d)}\|^2 \leq s a(n)}\Big) \bigg)
    \Bigg\|,
\end{align*}
and
\begin{align*}
  Z_2^{(d)} (s) := \max_{1\leq k\leq  n} \bigg\| \sum_{j =1}^{k}
  \E \Big(X_j^{(d)} \one_{\| X_j^{(d)} \|^2 \leq s a(n)}\Big)
  -\frac{k}{ n} \sum_{j=1}^{ n} \E \Big(X_j^{(d)} \one_{\|X_j^{(d)} \|^2 \leq s a(n)}\Big) \bigg\|.
\end{align*}
By the triangle inequality,
\begin{multline*}
  Z_1^{(d)} (s) \leq \max_{1\leq k\leq  n} \Bigg\| \sum_{j =1}^{k}
  \bigg( X_j^{(d)} \one_{\| X_j^{(d)} \|^2 \leq sa(n)}
  -\E \Big(X_j^{(d)} \one_{\| X_j^{(d)} \|^2 \leq sa(n)}\Big)\bigg)\Bigg\| \\
  + \Bigg\| \sum_{j=1}^{ n} \bigg(X_j^{(d)} \one_{\| X_j^{(d)}\|^2
    \leq s a(n)}
  -\E\Big(X_j^{(d)} \one_{\| X_j^{(d)} \|^2 \leq s a(n)}\Big) \bigg)\Bigg\|.
\end{multline*}
Note that
\begin{equation*}
  \Bigg\| \sum_{j=1}^{k} \bigg( X_j^{(d)} \one_{\| X_j^{(d)}\|^2\leq s a(n)}
  -\E \Big(X_j^{(d)} \one_{\| X_j^{(d)}\|^2 \leq s a(n)}\Big)\bigg) \Bigg\|^2,\quad k\in\N,
\end{equation*}
is a submartingale with respect to the filtration generated by the
sequence $\big(X_k^{(d)}\big)_{k \in \N}$. By Doob's inequality,
\begin{align*}
  \P\Big\{Z_1^{(d)} \geq 2^{-1} \eps \sqrt{a(n)}\Big\} 
  &\leq \Prob{\max_{1\leq k\leq  n} \Bigg\| \sum_{j=1}^{k}
      \bigg(X_j^{(d)} \one_{\| X_j^{(d)} \|^2 \leq s a(n)}
      -\E \Big(X_j^{(d)}\one_{\| X_j^{(d)} \|^2 \leq sa(n)}\Big)\bigg)\Bigg\|
    \geq 4^{-1} \eps \sqrt{a(n)}} \\
  & + \Prob{\Bigg\| \sum_{j=1}^{ n}\bigg( X_j^{(d)} \one_{\|X_j^{(d)} \|^2 \leq s a(n)}
      -\E\Big(X_j^{(d)} \one_{\| X_j^{(d)} \|^2 \leq sa(n)}\Big)\bigg)\Bigg\|
    \geq 4^{-1} \eps \sqrt{a(n)}} \\
  & \leq 32\eps^{-2}a(n)^{-1} \E\Bigg\| \sum_{j=1}^{ n}
      \bigg( X_j^{(d)} \one_{\| X_j^{(d)}\|^2 \leq s a(n)}-
      \E \Big(X_j^{(d)} \one_{\| X_j^{(d)} \|^2\leq s a(n)}\Big)
      \bigg)\Bigg\|^2 \\
  & =32\eps^{-2}a(n)^{-1}  n \E \bigg\|
      X_1^{(d)} \one_{\|X_1^{(d)} \|^2 \leq s a(n)}
      -\E \Big(X_1^{(d)} \one_{\| X_1^{(d)}\|^2 \leq s a(n)}\Big)\bigg\|^2 \\
  & \leq 64\eps^{-2}a(n)^{-1} n
  \E \Big(\| X_1^{(d)} \|^2 \one_{\| X_1^{(d)}\|^2 \leq s a(n)}\Big).
\end{align*}
Furthermore,
\begin{displaymath}
  Z_2^{(d)} \leq \max_{1\leq k\leq  n} \bigg\| \sum_{j=1}^{k}
  \E \Big(X_j^{(d)} \one_{\| X_j^{(d)} \|^2 \leq s a(n)} \Big)\bigg\|
  +\bigg\| \sum_{j=1}^{ n} \E
  \Big(X_j^{(d)} \one_{\| X_j^{(d)} \|^2\leq s a(n)} \Big) \bigg\| 
  \leq 2 n \bigg\| \E\Big(X_1^{(d)} \one_{\| X_1^{(d)} \|^2 \leq s a(n)}\Big)\bigg\|.
\end{displaymath}
Finally, assumptions \eqref{equ_ass_i_2} and \eqref{equ_ass_iii} imply
that
\begin{displaymath}
  \lim_{s \to 0+} \limsup_{d \to \infty} \Prob{Z_1^{(d)} \geq 2^{-1}\eps \sqrt{a(n)}} =0, 
  \quad\lim_{s \to 0+} \limsup_{d \to\infty}a(n)^{-1/2}Z_2^{(d)}(s) =0.
\end{displaymath}

\section*{Acknowledgements}
The author is grateful to Ilya Molchanov and Andrii Ilienko for their suggestions and
patience in instructing and correcting this paper, and the anonymous referee who 
draws the author's attention to the alternative continuous mapping theorem.

\bibliographystyle{cas-model2-names}



\end{document}